\documentclass{amsart}

\usepackage{amssymb}
\usepackage[pdftex]{graphicx}

\allowdisplaybreaks

\newtheorem{theorem}{Theorem}[section]
\newtheorem{lemma}[theorem]{Lemma}

\newtheorem{remark}[theorem]{Remark}

\newtheorem{TheoA}{Theorem A}
\newtheorem{TheoB}{Theorem B}

\theoremstyle{definition}

\newcommand{\N}{\mathbb{N}}
\newcommand{\Z}{\mathbb{Z}}
\newcommand{\R}{\mathbb{R}}
\newcommand{\C}{\mathbb{C}}
\newcommand{\X}{\mathcal{X}}

\newcommand{\A}{\mathcal{A}}

\newcommand{\CZ}{Cal\-de\-r\'on-Zyg\-mund }

\newcommand{\summ}{\sum\nolimits}

\newcommand{\dem}{\noindent {\bf Proof. }}
\newcommand{\demA}{\noindent {\bf Proof of Theorem A. }}
\newcommand{\demB}{\noindent {\bf Proof of Theorem B. }}
\newcommand{\fin}{\hspace*{\fill} $\square$ \vskip0.2cm}

\begin{document}

\title[Dyadic coverings and nondoubling CZ theory]
{A note on dyadic coverings and \\ nondoubling Calder\'on-Zygmund theory}

\author[Conde]
{Jose M. Conde}

\begin{abstract}
We construct a family of $n+1$ dyadic filtrations in $\R^n$, so that every Euclidean ball $B$ is contained in some cube $Q$ of our family satisfying $\mathrm{diam}(Q) \le c_n \mathrm{diam}(B)$ for some dimensional constant $c_n$. Our dyadic covering is optimal on the number of filtrations and improves previous results of Christ and Garnett/Jones by extending a construction of Mei for the $n$-torus. Based on this covering and motivated by applications to matrix-valued functions, we provide a dyadic nondoubling Calder\'on-Zygmund decomposition which avoids Besicovitch type coverings in Tolsa's decomposition. We also use a recent result of Hyt\"onen and Kairema to extend our dyadic nondoubling decomposition to the more general setting of upper doubling metric spaces.
\end{abstract}

\maketitle

\section*{Introduction}
\label{Section1}

\footnote{After this paper was written, Tuomas Hyt\"onen kindly communicated to us a different argument for Theorem A which he had included in the lecture notes of one of his courses.}

\addtolength{\parskip}{+1ex}

A dyadic system or filtration in $\R^n$ is a family $\bigcup_{k \in \Z} \mathcal{A}(k) \, = \, \bigcup_{k \in \Z} \, \bigcup_{j \in \Z^n} \{ Q_{jk} \}$ \\ [-0.6pt] where the $Q_{jk}$'s are (half open) cubes with sides parallel to the axes and satisfying the following properties:
\begin{itemize}
\item The side length $\ell(Q_{jk})$ of $Q_{jk}$ is $2^{-k}$,

\item The $k$-th generation $\mathcal{A}(k) = \bigcup_j \{ Q_{jk} \}$ is a partition of $\R^n$,

\item Any two cubes in the filtration are either disjoint or nested.
\end{itemize}
The starting point of this note is the following dyadic covering result.

\begin{TheoA}
There exist $n+1$ dyadic filtrations $\mathcal{A}_0, \mathcal{A}_1, \ldots, \mathcal{A}_n$ of $\R^n$ such that every Euclidean ball $B$ is contained in some cube $Q \in \bigcup_m \mathcal{A}_m$ satisfying that $\mathrm{diam}(Q) \le c_n \mathrm{diam}(B)$. The number of dyadic filtrations is optimal.
\end{TheoA}

Independently, Christ and Garnett/Jones provided the first dyadic covering of this kind by taking $3^n$ filtrations indexed by $s \in \{ 0, \frac13, \frac23 \}^n$ and with associated cubes determined by $$Q_{jk}^s = 2^{-k} \big( j + [0,1)^n + s \big).$$ Tao Mei reduced the number of filtrations to $n+1$ for the $n$-torus and also to 2 for the real line \cite{Mei,Mei2}. Nevertheless, his argument in \cite[Remark 7]{Mei} does not provide a dyadic covering of $\R^n$ for $n > 1$. Garnett/Jones \cite{GJ} combined dyadic systems to relate BMO with dyadic BMO. Theorem A shows that $\mathrm{BMO}_{\R^n}$ is isomorphic to the intersection of $n+1$ translates of dyadic BMO. Random dyadic systems have also been instrumental for Nazarov/Treil/Volberg's nonhomogeneous $T1$ theorem \cite{NTV2} or Hyt\"onen's solution of the $A_2$ conjecture \cite{Hytonen4}. 

Dyadic systems are crucial for a better understanding of many results in classical harmonic analysis. Geometrically, we obtain nested partitions of the space with unbeatable overlapping properties. Probabilistically, we may use the full machinery of $L_p$ martingale inequalities, usually exploited in the dyadic setting through the language of Haar functions. In some specific situations, the dyadic approach is the only one available. This is exactly the case for the matrix-valued Calder\'on-Zygmund theory, where the absence of Besicovitch or Vitali type coverings led to a dyadic approach introduced in \cite{Pa} through noncommutative martingales. Our work (in progress) on matrix-valued Calder\'on-Zygmund theory for nondoubling measures has motivated our Theorem B below, a dyadic form of Tolsa's nondoubling CZ decomposition \cite{Tolsa}. Namely, let $\mu$ be any measure on $\R^n$ with polynomial growth $\mu(B(x,r)) \le c_\mu r^k$ for some $0 < k \le n$. If $\mathcal{A} = \bigcup_m \mathcal{A}_m$ is the family of all the cubes in Theorem A, consider the maximal function $$\mathcal{M}_{\mathcal{A},\mu}f(x) \, = \, \sup_{x \in Q \in \mathcal{A}} \frac{1}{\mu(2Q)} \int_Q |f(y)| \, d\mu(y).$$ Then, we may write the $\lambda$-level set as an almost disjoint union $$\Omega_{\mathcal{A},\mu}^\lambda = \big\{ \mathcal{M}_{\mathcal{A},\mu}f > \lambda \big\} = \bigcup_{j \ge 1} Q_j,$$ with $Q_j \in \mathcal{A}$ maximal, i.e. there is no other $Q \in \mathcal{A}$ such that $Q_j \varsubsetneq Q \subset \Omega_{\mathcal{A},\mu}^\lambda$. The almost disjointness follows from maximality and the fact that $\mathcal{A}$ is a finite union of dyadic filtrations. Our dyadic nondoubling Calder\'on-Zygmund decomposition is partly formulated in terms of these cubes. Let us recall in passing that a cube $Q$ is called $(\alpha, \beta)$-doubling with respect to $\mu$ when $\mu(\alpha Q) \le \beta \mu(Q)$.

\begin{TheoB}
If $f \in L_1(\mu)$ and $\lambda > \|f\|_1/\|\mu\|$, let $\{Q_j\}$ be as above and define $w_j = \chi_{Q_j}/\sum_k \chi_{Q_k}$. Then, there exist $\alpha$ and $\beta$ ---depending only on $c_n$--- such that for every family of $\big( \alpha,\beta \big)$-doubling dyadic cubes $R_j \in \mathcal{A}$ with $3Q_j \subset R_j$, there exist functions $\varphi_j: \R^n \to \C$ satisfying the following properties$:$
\begin{itemize}
\item $\displaystyle \|\varphi_j\|_\infty \lesssim \frac{1}{\mu(R_j)} \int_{Q_j} f \, d\mu$ \, and \, $\displaystyle \Big\|\summ_j \varphi_j \Big\|_\infty \lesssim \lambda$.

\item $\displaystyle \mathrm{supp} \, \varphi_j \subset R_j$ \, and \, $\displaystyle \int \varphi_j \, d\mu = \int_{Q_j} f w_j \, d\mu \sim \int_{Q_j} f \, d\mu$.
\end{itemize}
A careful selection of $R_j$ yields a dyadic CZ decomposition adapted to $\mu$ $$f \, = \, g + b \, = \, \Big( f \chi_{\R^n \setminus \bigcup_j Q_j} + \summ_j \varphi_j \Big) + \summ_j \big( f w_j - \varphi_j \big).$$
\end{TheoB}

Our proof adapts Tolsa's construction in terms of our dyadic type covering in Theorem A. As in Tolsa's paper, we may prove the weak type $(1,1)$ boundedness of Calder\'on-Zygmund operators which satisfy standard $k$-dimensional size/smoothness kernel conditions and are bounded on $L_2(\mu)$. Moreover, we may also construct a dyadic form of Calder\'on-Zygmund decomposition in the (more general) setting of upper doubling metric spaces, which complements the approaches in \cite{Anh,Hytonen,NTV}. Our decomposition relies this time on Hyt\"{o}nen/Kairema's covering \cite{Hytonen2}. It is worth to mention that it would be quite interesting to refine such covering in terms of the number of selected filtrations, as we do in this paper for $\R^n$.  

%We end this section with some notes about notation: throughout the paper, the letter $C$ will denote a constant whose value may change every time it appears. The value of constants with subscripts, such as $C_K$, remains fixed. By a cube $Q$ we mean a cube with sides parallel to the axes. Its side length will be denoted by $l(Q)$, and its center by $x_Q$. In the same fashion, the center of a ball $B$ will be called $x_B$.

\section{The dyadic covering}
\label{Section2}

A dyadic filtration is completely determined by a distinguished cube (which we call initial cube) and its ancestors. Let us write $Q_0$ for the unit cube $[0,1)^n$ and $\mathbf{v}$ for the vector $(1,1,\ldots,1)$ in $\R^n$. Then, we construct precise dyadic filtrations $\mathcal{A}_0, \mathcal{A}_1, \ldots, \mathcal{A}_n$ for the proof of Theorem A as follows. 

\begin{enumerate}
\item Usual dyadic filtration $$\mathcal{A}_0 = \Big\{ Q_{jk}^0 = 2^{-k}(j+Q_0) \, \big| \ k \in \Z, \; j \in \Z^n \Big\}.$$

\item Shifted dyadic filtrations. Let $p_n$ be the smallest odd number strictly bigger than $n$. Then we set
$\mathcal{A}_m = \{ Q_{jk}^m \, | \ k \in \Z, \; j \in \Z^n\}$ for $1 \le m \le n$, where the cubes $Q_{jk}^m$ are determined by the following algorithm:

\vskip5pt

\begin{itemize}
\item The initial cube of $\mathcal{A}_m$ is $$Q_{00}^m= Q_0+ \frac{m}{p_n} \, \mathbf{v}.$$ This determines the k-th generation of $\mathcal{A}_m$ for any $k \ge 0$.

\vskip5pt

\item The parent of $Q_{0 0}^m$ is determined by 

\vskip3pt

\begin{itemize}
\item[i)] $Q_{0 -1}^m = Q_{0 -1}^0 + \lambda \mathbf{v}$

\item[ii)] The parameter $\lambda \in \frac{2}{p_n} \Z$.
\end{itemize}

\vskip3pt

Among the $2^n$ possible dyadic parents of $Q_{0 0}^m$, condition i) reduces our choice to only the two ones centered at the diagonal $\R \mathbf{v}$. Now, one and only one of these cubes fulfills condition ii) since $p_n$ is odd. This determines the $k$-th generation of $\mathcal{A}_m$ for $k=-1$. 

\vskip5pt

\item Inductively, the parent of $Q_{0k}^m$ is determined by $$Q_{0 \, k-1}^m = Q_{0\,k-1}^0 + \lambda \mathbf{v} \quad \mbox{for some} \quad \lambda \in 2^{-k+1}\frac{1}{p_n} \Z.$$ This determines the k-th generation of $\mathcal{A}_m$ for any $k < -1$.
\end{itemize}
\end{enumerate}

\demA With a slight change on the value of $c_n$, we may replace the ball $B$ by the circumscribed cube $R_B$ with sides parallel to the axes. Let us write $\ell(K)$ for the side length of the cube $K$. Fix $k_0 \in \Z$ determined by $$\frac{2^{-k_0-1}}{p_n} \le \ell(R_B) < \frac{2^{-k_0}}{p_n}.$$ We now claim that the set $$J_k = \frac{2^{-k}}{p_n} \Z^n$$ contains the vertices of all cubes in the $k$-th generation $\mathcal{A}_m(k)$ for all filtrations $0 \le m \le n$. Indeed, all vertices of $\mathcal{A}_m(0)$ are trivially in $J_0$ for any value of $m$. In addition, we observe that if a vertex of a cube in $\mathcal{A}_m(k)$ is in $J_k$, then all vertices of cubes in $\mathcal{A}_m(k)$ will be in $J_k$. If $k>0$ this is clear since the vertices of the initial cube of $\mathcal{A}_m$ are in $J_0 \subset J_k$. On the other hand, if $k<0$ the cube $Q_{0k}^m$ has a vertex at $\lambda \mathbf{v}$, which is in $J_k$ since $\lambda$ is an integer multiple of $2^{-k}/p_n$. This proves our claim. Now if $x,y \in J_{k_0}$ with $x \not= y$, we find $$|x-y|_\infty \ge \frac{2^{-k_0}}{p_n} > \ell(R_B).$$ Let us show that $R_B$ is contained in a cube of $\mathcal{A}_m(k_0)$ for some $0 \leq m \leq n$. Define $J_k(m) \subset J_k$ as the set of vertices of $\mathcal{A}_m(k)$ and let us write $P_j$ for the orthogonal projection onto the $j$-th axis of $\R^n$. We argue by contradiction, if $R_B$ were not contained in any cube of $\mathcal{A}_m(k_0)$, for each $0 \le m \le n$ we could find $1 \leq j_m \leq n$ such that $P_{j_m}(R_B) \cap P_{j_m}(J_{k_0}(m))$ is nonempty. By the pigeonhole principle, there would exist some $1 \le j \le n$ such that $P_{j}(R_B) \cap P_{j}(J_{k_0}(m))$ is nonempty for at least two values of $m$, say $m_1$ and $m_2$. Since $P_j(J_{k_0}(m_1))$ and $P_j(J_{k_0}(m_2))$ are disjoint, we obtain $\# (P_{j}(R_B) \cap P_{j}(J_{k_0})) \geq 2$. This contradicts the inequality $|x - y|_\infty > \ell(R_B)$ for $x,y \in J_{k_0}$ distinct and we conclude that $R_B \subset Q \in \mathcal{A}_m(k_0)$ for some $0 \le m \le n$. Moreover, the size of $Q$ is comparable to that of $R_B$ since $\mathrm{diam} (Q) = 2^{-k_0} \sqrt{n} \le 2 p_n \mathrm{diam}(R_B)$. Let us finally show that the number of dyadic filtrations is optimal. Given any family of $n$ dyadic filtrations $\mathcal{A}_1, \mathcal{A}_2, \ldots, \mathcal{A}_n$ of $\R^n$, we may always find $Q_m \in \mathcal{A}_m$ for $1 \le m \le n$ such that $$\bigcap_{m=1}^n \partial Q_m \neq \emptyset \quad \mbox{and} \quad \ell(Q_m) = 1.$$ Here $\partial Q_m$ denotes the boundary of $Q_m$. Taking $q \in \bigcap \partial Q_m$, the ball of radius $\varepsilon$ and center $q$ requires at least a cube of size length $2$. Let $\varepsilon \to 0$ to conclude. \fin 

%\begin{figure}[htb]
   %\centering
   %\includegraphics[scale=0.5]{diadicadesplazadacompleta.pdf}
   %\caption{\small{Initial cube and ancestors of $\mathcal{A}_0, \mathcal{A}_1, \mathcal{A}_2$ in dimension 2.}}
   %\label{picdiadicadesplazadacompleta}
%\end{figure}

%\begin{corollary}
%$\mathrm{BMO}_{\R^n}$ is isomorphic to the intersection of $n+1$ translates of the dyadic version of $\mathrm{BMO}_{\R^n}$, $\mathrm{BMO}_{\R^n}^d$.
%\end{corollary}
%See \cite{Mei} for the proof.

\begin{remark}
\emph{We say that $R_B$ \emph{fits} in $Q \in \A$ if $R_B \subset Q$ and $\mathrm{diam}(Q) \le c_n \mathrm{diam}(R_B)$.}
\end{remark}

\begin{remark}
\emph{Once we know how to construct dyadic ancestors of the initial cube for each filtration, we may produce many other families of $n+1$ dyadic filtrations satisfying Theorem A. We just need to follow Mei and shift the initial cubes making sure that their relative dyadic distances are positive, see \cite{Mei} for details.}
\end{remark}

\section{Dyadic nondoubling \\ Calder\'on-Zygmund decomposition}
\label{Section3}

We now use our dyadic covering to construct a dyadic form of Tolsa's nondoubling Calder\'on-Zygmund decomposition. Set $\Delta = \{(x,y) \in \R^{2n}: \; x=y \}$ and fix some $0 < k \le n$. A $k$-dimensional Calder\'on-Zygmund kernel is given by a function $K: \R^{2n} \setminus \Delta \to \C$ satisfying the standard size/smoothness conditions

\begin{itemize}
\item[a)] $\displaystyle |K(x,y)| \lesssim \frac{1}{|x-y|^k}$,

\vskip3pt

\item[b)] There exists $0 < \delta \le 1$ so that $|y-y'| \leq \frac{1}{2} |x-y|$ yields $$\big| K(x,y)-K(x,y') \big| \, + \, \big| K(y,x)-K(y',x) \big| \, \lesssim \, \frac{|y-y'|^\delta}{|x-y|^{k+\delta}}.$$
\end{itemize} 
We start adapting a couple of preliminary lemmas from \cite{Tolsa} to the dyadic setting. 

\begin{lemma} \label{sec2:lemma:doublingcubes}
Let $\mu$ be a measure on $\R^n$ of polynomial growth $\mu(B(x,r)) \le c_\mu r^k$ and let $\mathcal{A} = \bigcup_m \mathcal{A}_m$ and $c_n$ as in the statement of Theorem A. Then we may construct arbitrary small/large doubling dyadic cubes as follows\hskip1pt$:$
\begin{itemize}
\item[i)] Assume $\beta > (c_n \alpha)^n$. Then for $\mu$-almost every $x$ there is a sequence $\{Q_j\}$ of $(\alpha,\beta)$-doubling cubes with $x \in Q_j \in \A$ for all $j$ and such that $\ell(Q_j) < 2^{-j}$.

\vskip5pt

\item[ii)] Assume $\beta > (c_n \alpha)^k$. Then for $x \in \mathrm{supp}(\mu)$ and for each positive number $r$ there exists an $(\alpha,\beta)$-doubling cube with $x \in Q \in \mathcal{A}$ and such that $\ell(Q) > r$.
\end{itemize}
\end{lemma}

\dem It follows from Theorem A in conjunction with \cite[Remark 2.1 \& 2.2]{Tolsa}. \fin

\demB
Let us set $\alpha = 6 c_n$ and $\beta = (6 c_n^2)^k + 1$. We will assume with no loss of generality that $f$ is positive. We now adapt Tolsa's contruction to the dyadic setting as follows. Emulating the weak-$*$ limit argument in \cite{Tolsa}, we may suppose that the number of maximal dyadic cubes $\{Q_j\}$ for $\Omega_{\A,\mu}^\lambda$ is finite. Label the cubes $R_j$ so that their size lengths are nondecreasing. The functions $\varphi_j$ will have the form $\varphi_j= \gamma_j \chi_{A_j}$ with $$A_j \subset R_j \quad \mbox{and} \quad \gamma_j \in \R_+ \quad \mbox{so that} \quad \gamma_j \mu(A_j) = \int f w_j \, d\mu.$$ This provides $\mathrm{supp} \hskip1pt \varphi_j \subset R_j$ and $\int \varphi_j d\mu = \int_{Q_j} f w_j d\mu$. Now we claim that $$\frac{1}{\mu(A_j)} \, \le \, \frac{2}{\mu(R_j)} \qquad \mbox{and} \qquad \summ_j \varphi_j \, \le \, 4 \beta \lambda.$$ The remaining properties clearly follow from the claim. Taking $A_1 = R_1$ we find 
\begin{eqnarray*}
\varphi_1 \hskip7pt \le \hskip7pt \gamma_1 & = & \frac{1}{\mu(R_1)} \int f w_1 \, d\mu \\ & \le & \frac{\beta}{\mu(\alpha R_1)} \int_{R_1} f \, d\mu \hskip7pt \le \hskip7pt \frac{\beta}{\mu(2 R_1)} \int_{R_1} f \, d\mu \hskip7pt < \hskip7pt 4 \beta \lambda.
\end{eqnarray*}
Let us now proceed by induction. Assume that the family $\varphi_1, \varphi_2, \ldots, \varphi_{\ell-1}$ satisfies the claim above and consider the subfamily $\varphi_{s_1}, \varphi_{s_2}, \ldots, \varphi_{s_m}$ of functions with $R_{s_j} \cap R_\ell \neq \emptyset$. Since $\ell({R_{s_j}}) \le \ell(R_\ell)$ we have $R_{s_j} \subset 3R_\ell$ and         
\begin{eqnarray*}
\sum_{j=1}^m \int \varphi_{s_j} d\mu & = & \sum_{j=1}^m \int_{Q_{s_j}} f w_{s_j} d\mu \\ & \le & \int_{3R_\ell} fd\mu \hskip7pt \le \hskip7pt \lambda \mu(2S_\ell) \hskip7pt \le \hskip7pt \lambda \mu(\alpha R_\ell),
\end{eqnarray*}
where $S_\ell$ is a dyadic cube in $\A$ so that $3R_\ell$ fits in $S_\ell$. Chebychev's inequality gives $$\mu 
\Big\{ \summ_j \varphi_{s_j} > 2 \beta \lambda \Big\} \, \le \, \frac{\mu(\alpha R_\ell)}{2 \beta}.$$
Taking $A_\ell = R_\ell \cap \{ \sum_j \varphi_{s_j} \le 2 \beta \lambda \}$ we find 
\begin{eqnarray*}
\mu(A_\ell) & \ge & \mu(R_\ell) - \mu \Big\{ \sum_{j=1}^m \varphi_{s_j} > 2 \beta \lambda \Big\} \\ & \ge & \frac{\mu(\alpha R_\ell)}{\beta} - \frac{\mu(\alpha R_\ell)}{2 \beta} \hskip7pt \ge \hskip7pt \frac{\mu(R_\ell)}{2}. 
\end{eqnarray*}
On the other hand, it suffices to show that $\gamma_\ell \le 2 \beta \lambda$ since $$\sum_{j \le \ell} \varphi_j \, = \, \Big( \varphi_\ell + \sum_{j=1}^m \varphi_{s_j} \Big) \chi_{A_\ell} + \Big( \sum_{j \le \ell-1} \varphi_j \Big) \chi_{\R^n \setminus A_\ell} \, \le \, \big( \gamma_\ell + 2 \beta \lambda \big) \, \chi_{A_\ell} + 4 \beta \lambda \, \chi_{\R^n \setminus A_\ell}.$$ However, maximality of $Q_\ell$ gives $\displaystyle \gamma_\ell \hskip-1pt \le \hskip-1pt \frac{1}{\mu(A_\ell)} \int_{Q_\ell} f d\mu \le \frac{2 \beta}{\mu(\alpha R_\ell)} \int_{R_\ell} f d\mu \le 2 \beta \lambda$. \fin 

Given a $k$-dimensional CZ kernel $K$, a Calder\'on-Zygmund operator $T$ associated to it is any operator admitting the integral representation below for a suitable class of test functions $$Tf(x) = \int_{\R^n} K(x,y)f(y)d\mu(y) \quad \mbox{for} \quad x \not\in \mathrm{supp} f.$$ 
\begin{remark}
\emph{It is well-known that any CZO which is bounded on $L_2(\mu)$ satisfies the weak type $(1,1)$ estimate $T: L_1(\mu) \to L_{1,\infty}(\mu)$, see \cite{NTV,Tolsa}. Using the dyadic nondoubling CZ decomposition $f = g+b$ in Theorem B, we obtain a slightly different proof which appears to be the right approach for matrix-valued functions. Namely, this allows us to avoid Besicovitch type coverings, which do not seem to have easy noncommutative analogs. The proof follows from Theorem B as it followed in \cite{Tolsa} from the corresponding CZ decomposition, see the Appendix below for a more general setting. There is just one little deviation from Tolsa's argument here, which justifies in passing our choice $\beta = (6c_n^2)^k + 1$ in the proof of Theorem B, instead of the smaller $\beta = (6 c_n)^k + 1$. Indeed, let $Q$ and $R$ be two cubes in $\R^n$ such that $Q \subset R \in \mathcal{A}$ and there are no smaller $(\alpha,\beta)$-doubling cubes of $\A$ which contain $Q$. A key point in Tolsa's argument imposes here an upper bound for $$\int_{R\setminus Q} \frac{d\mu(x)}{|x-x_{Q}|^k}$$ independent of $Q$ and $R$, where $x_Q$ denotes the center of $Q$. This holds for $$\alpha>3 \quad \mbox{and} \quad \beta > (c_n \alpha)^k.$$ The proof is as follows. Let $Q_N = (c_n \alpha)^N Q$ be the smallest $(c_n \alpha, \beta)$-doubling cube of the form $Q_j = (c_n \alpha)^j Q$ for $j$ positive. Consider a cube $R' \in \A$ so that $Q_N$ fits in $R'$. Then $R'$ is $(\alpha,\beta)$-doubling since $\mu(\alpha R') \le \mu(c_n \alpha Q_N) \le \beta \mu(Q_N) \le \beta \mu(R')$ and contains $Q$. By hypothesis we have $R \subset R'$ and $$\ell(R) \le \ell(R') \le c_n \ell(Q_N).$$ Since $\alpha>3$ and $R \cap Q_N \not= \emptyset $, we find $R \subset 3 c_n Q_N \subset Q_{N+1}$ and $$\mu(Q_j) \le \frac{\mu(Q_{j+1})}{\beta} \le \cdots \le \frac{\mu(Q_N)}{\beta^{N-j}} \quad \mbox{for} \quad 0 \le j < N.$$
This, together with the assumption $\beta > (c_n \alpha)^k$ leads to the following estimate
\begin{eqnarray*}
\int_{R\setminus Q} \frac{d\mu(x)}{|x-x_{Q}|^k} & \leq & \sum_{j=1}^{N+1} \int_{Q_j \setminus Q_{j-1}} \frac{d\mu(x)}{|x-x_{Q}|^k} \, \lesssim \, \sum_{j=1}^{N+1} \frac{\mu(Q_j)}{\ell(Q_j)^k} \\ & \le &  \sum_{j=1}^N \, \frac{\beta^{j-N} \mu(Q_N)}{(c_n \alpha)^{(j-N)k} \ell(Q_N)^k} \ + \ \frac{\beta \mu(Q_{N})}{\ell(Q_{N+1})^k} \\ & \le & \frac{\mu(Q_N)}{\ell(Q_N)^k} \sum_{j=0}^\infty \frac{(c_n\alpha)^{jk} }{\beta^j} \ + \ \frac{\beta \mu(Q_{N})}{(c_n\alpha)^k \ell(Q_{N})^k} \ \lesssim \ \frac{\mu(Q_N)}{\ell(Q_N)^k}.
\end{eqnarray*}
Now the conclusion follows from the polynomial growth of the measure $\mu$. \fin} 
\end{remark}

%\appendix
\section*{Appendix. Metric spaces}
\label{Appendix}

\newtheorem{Atma}{Theorem}
\renewcommand{\theAtma}{A\arabic{Atma}}

A metric space $(\X,d)$ is called geometrically doubling if there exists some natural number $N \in \N$ such that for every $x \in \X$ and for every $r>0$, we may find $x_1, x_2, \ldots, x_N$ so that $B(x,r)$ can be covered by the balls $B(x_j,r/2)$. This holds iff for every $\delta \in (0,1)$, any ball $B(x,r) \subset\X$ can contain at most $N\delta^{-n}$ centres $z_j$ of disjoint balls $B(z_j,\delta r)$. The number $n = \log_2 N$ is the so-called geometric dimension of the space $\X$. We will assume that we have a Borel measure $\mu$ defined on $\X$ which is non-atomic and finite on bounded sets. The space $(\X,d,\mu)$ is said to be \emph{upper doubling} if there exists a function $\lambda: \X \times \R_+ \rightarrow \R_+$ such that $r \mapsto \lambda(x,r)$ is nondecreasing and verifies $$\mu\big( B(x,r) \big) \, \le \, \lambda(x,r) \, \le \, C_\lambda \, \lambda\left(x,\frac{r}{2}\right).$$ 
Both doubling measures and measures of polynomial growth are upper doubling. 

The analog of Christ and Garnett/Jones coverings has been recently found by Hyt\"onen and Kairema in \cite{Hytonen2}. Let us summarize their main result in the context of metric spaces; quasi-metric spaces can also be considered but we avoid it for simplicity of the exposition. Given an $n$-dimensional geometrically doubling space $(\X, d)$ and $\delta > 0$ small enough, we may find points $z_{jk}$ ---with $j$ running along certain index set $\Sigma_k$ for every $k \in \Z$--- and sets $Q_{jk}$ satisfying$\, :$  
\begin{enumerate}
\item $B(z_{jk}, \delta^k) \subset Q_{jk} \subset 48 B(z_{jk}, \delta^k)$,

\vskip3pt

\item $\bigcup_j Q_{jk}$ is a partition of $\mathcal{X}$ for all $k \in \Z$, 

\vskip3pt

\item If $k_1 \ge k_2 $, then either $Q_{j_1 k_1} \subset Q_{j_2 k_2}$ or $Q_{j_1 k_1} \cap Q_{j_2 k_2} = \emptyset$.
\end{enumerate}
Moreover, there exist a family of filtrations $$\A_m = \bigcup_{k \in \Z} \bigcup_{j \in \Sigma_k} Q_{jk}^m \quad \mbox{with} \quad 1 \le m \le K(n,\delta)$$ satisfying the properties above, such that every $d$-ball $B$ in $\mathcal{X}$ is contained in some $Q \in \A = \bigcup_m \A_m$ satisfying $\mathrm{diam}(Q) \le \frac12 C_\delta \, \mathrm{diam}(B)$. The first analogues of dyadic cubes in noneuclidean metric spaces go back to G. David and M. Christ. According to the terminology above, we say that $B$ fits in $Q \in \A$ when $B \subset Q \subset C_\delta B$. We have two comments:

\textbf{A.} It is interesting to note that the only constant above which depends on the geometric dimension $n$ is the number of filtrations, in contrast with the embedding constants 48 or $C_\delta$. A quick look at the argument in \cite{Hytonen2} shows that $K(n,\delta)$ is majorized by $(4/\delta)^{3n}$, which seems far to be optimal in view of Theorem A above. That leads us to propose the following

\vskip3pt

\noindent \textbf{Problem A.1.} Do we have dyadic coverings on $(\mathcal{X},d)$ with $K(n,\delta) \sim C_\delta n$?

\vskip3pt

\textbf{B.} Now we have all the tools to generalize Theorem B to the context of upper doubling metric spaces. We say that $Q \in \A$ is $(\alpha,\beta)$-doubling if there exists an $(\alpha,\beta)$-doubling $d$-ball $B$ fitting in $Q$. This notion could have also been used in the Euclidean case considered in the previous section. Then, Lemma \ref{sec2:lemma:doublingcubes} still holds in this setting since we know from \cite{Hytonen3} that arbitrary small/large $(\alpha,\beta)$-doubling concentric cubes can be constructed under suitable assumptions on the values of $\alpha, \beta$. We also consider a slightly different maximal function $$\mathcal{M}_{\mathcal{A},\mu}f(x) \, = \, \sup_{x \in Q \in \mathcal{A}} \, \sup_{B \subset Q \subset C_\delta B} \, \frac{1}{\mu(2 C_\delta B)} \int_Q |f(y)| \, d\mu(y).$$ We still have $\{ \mathcal{M}_{\mathcal{A},\mu}f > \lambda \} = \bigcup_{j \ge 1} Q_j$ with $Q_j \in \mathcal{A}$ maximal. According to \cite{Hytonen3}, a Calder\'on-Zygmund kernel in an upper doubling metric space $(\mathcal{X},d,\mu,\lambda)$ must satisfy

\begin{itemize}
\item[a)] $\displaystyle |K(x,y)| \lesssim \frac{1}{\lambda(x,d(x,y))}$,

\vskip3pt

\item[b)] There exists $0 < \delta \le 1$ so that $d(y,y') \leq \frac{1}{2} d(x,y)$ yields $$\hskip25pt \big| K(x,y)-K(x,y') \big| \, + \, \big| K(y,x)-K(y',x) \big| \, \lesssim \, \frac{d(y,y')^\delta}{d(x,y)^{\delta} \lambda(x,d(x,y))}.$$
\end{itemize} 

The result below generalizes Theorem B to the context of upper metric spaces and complements the approaches given in \cite{Anh,Hytonen,NTV}. The dyadic covering from \cite{Hytonen2} is essentially all what is needed to reproduce our argument in Section 2, thus we omit the details of the proof since they would be quite repetitive.  

\noindent \textbf{Theorem A.2.} \emph{Let $(\mathcal{X},d,\mu)$ be an upper doubling metric space. If $f \in L_1(\mathcal{X},\mu)$ and $\lambda > \|f\|_1/\|\mu\|$, let $w_j = \chi_{Q_j}/\sum_k \chi_{Q_k}$ with $\{Q_j\}$ as above. Consider a family of balls $B_j$ fitting in $Q_j$. Then, there exist $\alpha$ and $\beta$ ---depending only on $C_\delta$--- such that for every family of $\big( \alpha,\beta \big)$-doubling dyadic cubes $R_j \in \mathcal{A}$ with $3C_\delta B_j \subset R_j$ there exist functions $\varphi_j: \R^n \to \C$ satisfying$:$}
\begin{itemize}
\item \emph{$\displaystyle \|\varphi_j\|_\infty \lesssim \frac{1}{\mu(R_j)} \int_{Q_j} f \, d\mu$ \, and \, $\displaystyle \Big\|\summ_j \varphi_j \Big\|_\infty \lesssim \lambda$.}

\item \emph{$\displaystyle \mathrm{supp} \, \varphi_j \subset R_j$ \, and \, $\displaystyle \int \varphi_j \, d\mu = \int_{Q_j} f w_j \, d\mu \sim \int_{Q_j} f \, d\mu$.}
\end{itemize}
\emph{In particular, the CZ decomposition adapted to $(\mathcal{X},d,\mu)$ $$f \, = \, g + b \, = \, \Big( f \chi_{\mathcal{X} \setminus \bigcup_j Q_j} + \summ_j \varphi_j \Big) + \summ_j \big( f w_j - \varphi_j \big)$$ allows to show that any $L_2$-bounded CZO on $(\mathcal{X},d,\mu)$ maps $L_1(\mathcal{X},\mu)$ to $L_{1,\infty}(\mathcal{X},\mu)$.}

\vskip3pt

\noindent \textbf{Acknowledgement.} Supported in part by the ERC Grant StG-256997-CZOSQP (European Union) and also by the MEC Grant MTM-2010-16518 (Spain).  The author is grateful to Jos\'e Garc\'ia-Cuerva and Javier Parcet for useful conversations related to the content of this paper.

\bibliographystyle{amsplain}

\

\hfill \noindent \textbf{Jose M. Conde} \\
\null \hfill Instituto de Ciencias Matem\'aticas \\ 
\null \hfill CSIC-UAM-UC3M-UCM \\
\null \hfill Consejo Superior de Investigaciones Cient\'ificas \\
\null \hfill C/ Nicol\'as Cabrera 13-15. 28049, Madrid. Spain \\
\null \hfill\texttt{jose.conde@icmat.es}

\end{document}